# $E$-Orlicz Spaces and $E$-Orlicz-Sobolev space

**Abdulhameed Qahtan Abbood Altai**[1]

College of Islamic Sciences

University of Babylon, 51002, Babil, Iraq

and

**Nada Mohammed Abbas Alsafar**[2]

College of Education for Pure Sciences

University of Babylon, 51002, Babil, Iraq

**Abstract.** A new class of $E$-convex functions called $E$-$N$-functions, $E$-Young functions, $E$-strong Young functions and $E$-Orlicz functions are introduced by relaxing the definitions of $N$-functions, Young functions, strong Young functions and Orlicz functions. Then, new classes of $E$-Orlicz spaces and $E$-Orlicz-Sobolev space will be defined.

**1. Introduction.** Youness introduced the concept of $E$-convex sets and $E$-convex functions to generalize the concept of convex sets and convex functions to extend the studying of the optimality for non linear programming problems in 1999 [6]. Chen defined the semi-$E$-convex functions and studied its basic properties in 2002 [2]. Syau and Lee discussed properties of $E$-convex functions and introduced the concepts of pseudo $E$-convex functions and $E$-quasiconvex functions and strictly $E$-quasiconvex functions in 2004 [5]. The concept of Semi strongly $E$-convex functions was introduced by Youness and Tarek Emam in 2005 [7]. Sheiba Grace and Thangavelu considered the algebraic properties of $E$-convex sets in 2009 [3]. $E$-differentiable convex functions was defined by Meghed, Gomma, Youness and El-Banna [4] to transform a non-differentiable function to a differentiable function in 2013. Semi-$E$-Convex Function were introduced in 2015 by Ayache and Khaled [1]. Our contribution in this paper is to define the $E$-$N$-functions, $E$-Young functions, $E$-strong Young functions and $E$-Orlicz functions using the concept of $E$-convex functions to transfer the functions that are not in the class of $N$-functions, Young functions, strong Young functions and Orlicz functions into class of $N$-functions, Young functions, strong Young functions and Orlicz functions respectively. Then, the $E$-Orlicz spaces and $E$-Orlicz-Sobolev spaces will be defined using these new concepts.

**2. Definitions of $E$-convex functions $\Phi$**

**Definition 2.1**[6]: A set $S \subset R^n$ is said to be $E$-convex iff there is a map $E: R^n \to R^n$ such that $\lambda E(x) + (1-\lambda) E(y) \in S$, for each $x, y \in S, 0 \le \lambda \le 1$.

**Definition 2.2**[6]: A function $f: R^n \to R$ is said to be $E$-convex on a set $S \subset R^n$ iff there is a map $E: R^n \to R^n$ such that $S$ is an $E$-convex set and

$$f(\lambda E(x) + (1-\lambda)E(y)) \le \lambda f(E(x)) + (1-\lambda) f(E(y)),$$

[1]E-mail: ahbabil1983@gmail.com,

[2]E-mail: nadaalsafar333@gmail.com

$\forall x, y \in M, 0 \leq \lambda \leq 1$. And $f$ is called $E$-concave on a set $S$ if

$$f(\lambda E(x) + (1-\lambda)E(y)) \geq \lambda f(E(x)) + (1-\lambda)f(E(y)),$$

$\forall x, y \in S, 0 \leq \lambda \leq 1$.

**Definition 2.3**: Let $(\Omega, \Sigma, \mu)$ be a measure space. A function $\Phi: \Omega \times [0, \infty) \to \mathbb{R}$ is called an $E$-$N$-function if there exists a map $E: \Omega \times [0, \infty) \to \Omega \times [0, \infty)$ such that for $\mu$-a.e. $t \in \Omega$, $[0, \infty)$ is $E$-convex, $\Phi(E(t, u))$ is even continuous convex of $u$ on $[0, \infty)$, $\Phi(E(t, u)) > 0$ for any $u \in (0, \infty)$, $\lim_{u \to 0^+} \frac{\Phi(E(t,u))}{u} = 0$, $\lim_{u \to \infty} \frac{\Phi(E(t,u))}{u} = \infty$ and for each $u \in [0, \infty), \Phi(E(t, u))$ is a $\mu$-measurable function of $t$ on $\Omega$.

**Remark 2.1**: Every $N$-function is $E$-$N$-function where $E$ is the identity map.

**Examples 2.1**:

1. Let $\Phi: \mathbb{R} \times [0, \infty) \to \mathbb{R}$ be defined as $\Phi(t, u) = tu^2$ and let $E: \mathbb{R} \times [0, \infty) \to \mathbb{R} \times [0, \infty)$ be defined as $E(t, u) = (|t|, u)$. Then $\Phi$ is $E$-$N$-function but is not $N$-function because that for $\mu$-a.e. $t \in \mathbb{R}, \Phi(t, u)$ is concave of $u$ for $t \in (-\infty, 0)$.
2. Let $\Phi: \mathbb{R} \times [0, \infty) \to \mathbb{R}$ be defined as $\Phi(t, u) = (1-t)u^2 + t\exp(u)$ and let $E: \mathbb{R} \times [0, \infty) \to \mathbb{R} \times [0, \infty)$ be defined as $E(t, u) = (t, \ln u^2)$. So, $\Phi$ is $E$-$N$-function but is not $N$-function since for $\mu$-a.e. $t \in \mathbb{R}, \Phi(t, u)$ is not even.

**Definition 2.4**: Let $(\Omega, \Sigma, \mu)$ be a measure space. A function $\Phi: \Omega \times [0, \infty) \to \mathbb{R}$ is called an $E$-Young function if there exists a map $E: \Omega \times [0, \infty) \to \Omega \times [0, \infty]$ such that for $\mu$-a.e. $t \in \Omega$, $[0, \infty)$ is $E$-convex, $\Phi(E(t, u))$ is convex of $u$ on $[0, \infty)$, $\Phi(E(t, 0)) = \lim_{u \to 0^+} \Phi(E(t, u)) = 0$, $\lim_{u \to \infty} \Phi(E(t, u)) = \infty$ and for each $u \in [0, \infty)$, $\Phi(E(t, u))$ is a $\mu$-measurable function of $t$ on $\Omega$.

**Remark 2.2**: Every Young function is $E$-Young function where $E$ is the identity map.

**Examples 2.2**:

1. Let $\Phi: \mathbb{R} \times [0, \infty) \to \mathbb{R}$ be defined as $\Phi(t, u) = e^{t+u} - 1$ and let $E: \mathbb{R} \times [0, \infty) \to \mathbb{R} \times [0, \infty)$ be a map defined as $E(t, u) = (u, u)$. Then, $\Phi$ is $E$-Young function but is not Young function because for $\mu$-a.e. $t \in \mathbb{R}, \Phi(t, 0) = e^t - 1 \neq 0$.

2. Let $\Phi: \mathbb{C} \times [0, \infty) \to \mathbb{R}$ be defined as
$$\Phi(t, u) = \begin{cases} t\ln(u), & u > 1 \\ 0, & 0 \leq u \leq 1 \end{cases}$$
and let $E: \mathbb{C} \times [0, \infty) \to \mathbb{C} \times [0, \infty)$ be a map defined as $E(t, u) = (-|t|, u)$. So, $\Phi$ is $E$-Young function but is not Young function because that for $\mu$-a.e. $t \in \mathbb{C}, \Phi(t, u)$ is not convex because for $t \in (0, \infty)$ that $\frac{\partial^2 \Phi}{\partial u^2} = -\frac{t}{u^2} < 0$.

**Definition 2.5**: Let $(\Omega, \Sigma, \mu)$ be a measure space. A function $\Phi: \Omega \times [0, \infty) \to \mathbb{R}$ is called an $E$-strong Young function if there exists a map $E: \Omega \times [0, \infty) \to \Omega \times [0, \infty)$ such that for $\mu$-a.e. $t \in \Omega$, $[0, \infty)$ is $E$-convex, $\Phi(E(t, u))$ is convex continuous of $u$ on $[0, \infty)$, $\Phi(E(t, 0)) = 0 \Leftrightarrow u = 0$, $\lim_{u \to \infty} \Phi(E(t, u)) = \infty$ and for each $u \in [0, \infty), \Phi(E(t, u))$ is a $\mu$-measurable function of $t$ on $\Omega$.

**Remark 2.3**: Every strong Young function is $E$-strongYoung function where $E$ is the identity map.

**Examples 2.3**:

1. Let $\Phi: \mathbb{R} \times [0, \infty) \to \mathbb{R}$ be defined as $\Phi(t, u) = e^{u^t} - 1$ and let $E: \mathbb{R} \times [0, \infty) \to \mathbb{R} \times [0, \infty)$ be defined as $E(t, u) = (|t|, u)$. Then $\Phi$ is $E$-strong Young function but is not strong Young function because that $\Phi(t, u) = e^{u^t} - 1$ is not convex because for $t \in (-\infty, 0)$ that $u^t$ is not convex.
2. Let $\Phi: [0, \infty) \times [0, \infty) \to \mathbb{R}$ be defined as $\Phi(t, u) = \cosh(te^u) - 1$ and let $E: [0, \infty) \times [0, \infty) \to [0, \infty) \times [0, \infty)$ be defined as $E(t, u) = (u, 0)$. Then $\Phi$ is $E$-strong Young function but is not strong Young function because for $\mu$-a.e. $t \in [0, \infty)$, $\Phi(t, 0) = \cosh(t) - 1 \neq 0$.

**Definition 2.6**: Let $(\Omega, \Sigma, \mu)$ be a measure space. A function $\Phi: \Omega \times [0, \infty) \to \mathbb{R}$ is called an $E$-Orlicz function if there exists a map $E: \Omega \times [0, \infty) \to \Omega \times [0, \infty)$ such that for $\mu$-a.e. $t \in \Omega$, $[0, \infty)$ is $E$-convex, $\Phi(E(t, u))$ is convex of $u$ on $[0, \infty)$, $\Phi(E(t, 0)) = 0$, $\lim\limits_{u \to \infty} \Phi(E(t, u)) = \infty$, $0 < \Phi(E(t, u)) < \infty$ for any $u \in (0, \infty)$ and $\Phi(E(t, u))$ is left continuous at

$$U_\Phi = \sup\{u > 0: \Phi(E(t, u)) < +\infty\}.$$

Furthermore, for each $u \in [0, \infty)$, $\Phi(E(t, u))$ is a $\mu$-measurable function of $t$ on $\Omega$.

**Remark 2.4**: Every Orlicz function is $E$-Orlicz function where $E$ is the identity map.

**Examples 2.4**:

1. Let $\Phi: \mathbb{R} \times [0, \infty) \to \mathbb{R}$ be defined as $\Phi(t, u) = -t + u$ and let $E: \mathbb{R} \times [0, \infty) \to \mathbb{R} \times [0, \infty)$ be defined as $E(t, u) = (0, u^p), p \geq 1$. Then $\Phi$ is $E$-Orlicz function but is not Orlicz function. because for $\mu$-a.e. $t \in \mathbb{R}, \Phi(t, 0) = -t \neq 0$.
2. Let $\Phi: \mathbb{R} \times [0, \infty) \to \mathbb{R}$ be defined as $\Phi(t, u) = t + u^{p/(1-t)}, p \geq 1$ and let $E: \mathbb{R} \times [0, \infty) \to \mathbb{R} \times [0, \infty)$ be defined as $E(t, u) = (0, u)$. Then $\Phi$ is $E$-Orlicz function but is not Orlicz function because for $\mu$-a.e. $t \in \mathbb{R}, \Phi(t, 0) = t \neq 0$.

### 3. Elementary Properties of $E$-convex functions $\Phi$

**Proposition 3.1.** Let $(\Omega, \Sigma, \mu)$ be a measure space and $\Phi_1, \Phi_2: \Omega \times [0, \infty) \to \mathbb{R}$ be $E$-N-functions with a map $E: \Omega \times [0, \infty) \to \Omega \times [0, \infty)$. Then $\Phi_1 + \Phi_2$ and $c\Phi_1, c \geq 0$ are $E$-N-functions with $E$.

**Proposition 3.2.** Let $(\Omega, \Sigma, \mu)$ be a measure space with $\Omega \times [0, \infty)$ compact. If $\Phi_n: \Omega \times [0, \infty) \to \mathbb{R}$ is a sequence of continuous $E$-N-functions with a map $E: \Omega \times [0, \infty) \to \Omega \times [0, \infty)$ and $\Phi_n$ converges uniformaly to a continuous function $\Phi: \Omega \times [0, \infty) \to \mathbb{R}$, then $\Phi$ is an $E$-N-function with $E$.

**Proof.** Since $\Phi_n$ are continuous $E$-N-functions with $E$ such that $\Phi_n \to \Phi$ uniformly on compact set $\Omega \times [0, \infty)$ and $\Phi$ is continuous on $\Omega \times [0, \infty)$, then $\Phi_n(E) \to \Phi(E)$ uniformaly on $\Omega \times [0, \infty)$. So, for $\mu$-a.e. $t \in \Omega, \Phi(E(t, u)) = \lim\limits_{n \to \infty} \Phi_n(E(t, u))$ is even continuous convex of $u$ on $[0, \infty), \Phi(E(t, u)) > 0$ for any $u \in (0, \infty)$,

$$\lim_{u \to 0} \frac{\Phi(E(t,u))}{u} = \lim_{n \to \infty} \lim_{u \to 0} \frac{\Phi_n(E(t,u))}{u} = 0,$$

$$\lim_{u \to \infty} \frac{\Phi(E(t,u))}{u} = \lim_{n \to \infty} \lim_{u \to \infty} \frac{\Phi_n(E(t,u))}{u} = \infty$$

and for each $u \in [0, \infty)$, $\Phi(E(t,u))$ is a $\mu$-measurable function of $t$ on $\Omega$. ∎

**Proposition 3.3.** Let $(\Omega, \Sigma, \mu)$ be a measure space and $\Phi: \Omega \times [0, \infty) \to \mathbb{R}$ be a linear $E$-$N$-function with maps $E_1, E_2: \Omega \times [0, \infty) \to \Omega \times [0, \infty)$. Then $\Phi$ is an $E$-$N$-function with $E_1 + E_2$ and $cE_1$, $c \geq 0$ respectively.

**Proposition 3.4.** Let $(\Omega, \Sigma, \mu)$ be a measure space with $\Omega \times [0, \infty)$ compact. If $\Phi: \Omega \times [0, \infty) \to \mathbb{R}$ is a continuous $E$-$N$-function with a sequence of maps $E_n: \Omega \times [0, \infty) \to \Omega \times [0, \infty)$ that converges uniformaly to a map $E: \Omega \times [0, \infty) \to \Omega \times [0, \infty)$, then $\Phi$ is an $E$-$N$-function with $E$.

**Proof.** Because $\Phi$ is a continuous $E$-$N$-function with $E_n$ such that $E_n \to E$ uniformly on a compact set $\Omega \times [0, \infty)$, then $\Phi(E_n) \to \Phi(E)$ uniformaly on $\Omega \times [0, \infty)$. So, for $\mu$-a.e. $t \in \Omega$, $\Phi(E(t,u)) = \lim_{n \to \infty} \Phi(E_n(t,u))$ is even continuous convex of $u$ on $[0, \infty)$, $\Phi(E(t,u)) > 0$, $u \in (0, \infty)$,

$$\lim_{u \to 0} \frac{\Phi(E(t,u))}{u} = \lim_{n \to \infty} \lim_{u \to 0} \frac{\Phi(E_n(t,u))}{u} = 0,$$

$$\lim_{u \to \infty} \frac{\Phi(E(t,u))}{u} = \lim_{n \to \infty} \lim_{u \to \infty} \frac{\Phi(E_n(t,u))}{u} = \infty$$

and for each $u \in [0, \infty)$, $\Phi(E(t,u))$ is a $\mu$-measurable function of $t$ on $\Omega$. ∎

**Proposition 3.5.** Let $(\Omega, \Sigma, \mu)$ be a measure space with $\Omega \times [0, \infty)$ compact. If $\Phi_n: \Omega \times [0, \infty) \to \mathbb{R}$ is a sequence of continuous $E$-$N$-functions with a sequence of continuous maps $E_n: \Omega \times [0, \infty) \to \Omega \times [0, \infty)$ such that $\Phi_n$ converges uniformaly to a continuous function $\Phi: \Omega \times [0, \infty) \to \mathbb{R}$ and $E_n$ converges uniformaly to a continuous map $E: \Omega \times [0, \infty) \to \Omega \times [0, \infty)$, then $\Phi$ is an $E$-$N$-function with $E$.

**Proof.** We have $\Phi_n$ are continuous $E$-$N$-functions with continuous maps $E_n$ such that $\Phi_n \to \Phi$ and $E_n \to E$ uniformly on a compact set $\Omega \times [0, \infty)$ and $\Phi$ and $E$ are continuous on $\Omega \times [0, \infty)$. So, $\Phi_n(E_n) \to \Phi(E)$ uniformaly on $\Omega \times [0, \infty)$. Then, for $\mu$-a.e. $t \in \Omega$, that $\Phi(E(t,u)) = \lim_{n \to \infty} \Phi_n(E_n(t,u))$ is even continuous convex of $u$ on $[0, \infty)$, $\Phi(E(t,u)) > 0$, $u \in (0, \infty)$,

$$\lim_{u \to 0} \frac{\Phi(E(t,u))}{u} = \lim_{n \to \infty} \lim_{u \to 0} \frac{\Phi_n(E_n(t,u))}{u} = 0,$$

$$\lim_{u \to \infty} \frac{\Phi(E(t,u))}{u} = \lim_{n \to \infty} \lim_{u \to \infty} \frac{\Phi_n(E_n(t,u))}{u} = \infty$$

and for each $u \in [0, \infty)$, $\Phi(E(t,u))$ is a $\mu$-measurable function of $t$ on $\Omega$. ∎

**Proposition 3.6.** Let $(\Omega, \Sigma, \mu)$ be a measure space and $\Phi_1, \Phi_2: \Omega \times [0, \infty) \to \mathbb{R}$ are $E$-Young functions with a map $E: \Omega \times [0, \infty) \to \Omega \times [0, \infty)$. Then both $\Phi_1 + \Phi_2$ and $c\Phi_1, c \geq 0$ are $E$-Young functions with $E$.

**Proposition 3.7.** Let $(\Omega, \Sigma, \mu)$ be a measure space with $\Omega \times [0, \infty)$ compact. If $\Phi_n: \Omega \times [0, \infty) \to \mathbb{R}$ is a sequence of continuous $E$-Young functions with a map $E: \Omega \times [0, \infty) \to \Omega \times [0, \infty)$ and $\Phi_n$ converges uniformaly to a continuous map $\Phi: \Omega \times [0, \infty) \to \mathbb{R}$, then $\Phi$ is an $E$-Young function with $E$.

**Proof.** Since $\Phi_n$ are continuous $E$-Young functions with $E$ such that $\Phi_n \to \Phi$ uniformly on a compact set $\Omega \times [0, \infty)$ and $\Phi$ is continuous on $\Omega \times [0, \infty)$, then $\Phi_n(E) \to \Phi(E)$ uniformaly on $\Omega \times [0, \infty)$. Therefore, for $\mu$-a.e. $t \in \Omega$, $\Phi(E(t, u)) = \lim_{n \to \infty} \Phi_n(E(t, u))$ is convex of $u$ on $[0, \infty)$,
$$\Phi(E(t, 0)) = \lim_{u \to 0^+} \Phi(E(t, u)) = \lim_{n \to \infty} \lim_{u \to 0^+} \Phi_n(E(t, u)) = 0,$$
$$\lim_{u \to \infty} \Phi(E(t, u)) = \lim_{n \to \infty} \lim_{u \to \infty} \Phi_n(E(t, u)) = \infty$$
and for each $u \in [0, \infty)$, $\Phi(E(t, u))$ is a $\mu$-measurable function of $t$ on $\Omega$. ∎

**Proposition 3.8.** Let $(\Omega, \Sigma, \mu)$ be a measure space and $\Phi: \Omega \times [0, \infty) \to \mathbb{R}$ be a linear $E$-Young function with maps $E_1, E_2: \Omega \times [0, \infty) \to \Omega \times [0, \infty)$. Then $\Phi$ is $E$-Young functions with $E_1 + E_2$ and $cE_1, c \geq 0$ respectively.

**Proposition 3.9.** Let $(\Omega, \Sigma, \mu)$ be a measure space with $\Omega \times [0, \infty)$ compact. If $\Phi: \Omega \times [0, \infty) \to \mathbb{R}$ is a continuous $E$-Young function with a sequence of maps $E_n: \Omega \times [0, \infty) \to \Omega \times [0, \infty)$ that converges uniformaly to a map $E: \Omega \times [0, \infty) \to \Omega \times [0, \infty)$, then $\Phi$ is an $E$-Young function with $E$.

**Proof.** Because $\Phi$ is a continuous $E$-Young function with $E_n$ such that $E_n \to E$ uniformaly on a compact set $\Omega \times [0, \infty)$, then $\Phi(E_n) \to \Phi(E)$ uniformaly on $\Omega \times [0, \infty)$. Thus, for $\mu$-a.e. $t \in \Omega$, $\Phi(E(t, u)) = \lim_{n \to \infty} \Phi(E_n(t, u))$ is convex of $u$ on $[0, \infty)$,
$$\Phi(E(t, 0)) = \lim_{u \to 0^+} \Phi(E(t, u)) = \lim_{n \to \infty} \lim_{u \to 0^+} \Phi(E_n(t, u)) = 0,$$
$$\lim_{u \to \infty} \Phi(E(t, u)) = \lim_{n \to \infty} \lim_{u \to \infty} \Phi(E_n(t, u)) = \infty$$
and for each $u \in [0, \infty)$, $\Phi(E(t, u))$ is a $\mu$-measurable function of $t$ on $\Omega$. ∎

**Proposition 3.10.** Let $(\Omega, \Sigma, \mu)$ be a measure space with $\Omega \times [0, \infty)$ compact. If $\Phi_n: \Omega \times [0, \infty) \to \mathbb{R}$ is a sequence of continuous $E$-Young functions with a sequence of continuous maps $E_n: \Omega \times [0, \infty) \to \Omega \times [0, \infty)$ such that $\Phi_n$ converges uniformaly to a continuous function $\Phi: \Omega \times [0, \infty) \to \mathbb{R}$ and $E_n$ converges uniformaly to a continuous map $E: \Omega \times [0, \infty) \to \Omega \times [0, \infty)$, then $\Phi$ is an $E$-Young function with $E$.

**Proof.** We have $\Phi_n$ are continuous $E$-Young functions with continuous maps $E_n$ such that $\Phi_n \to \Phi$ and $E_n \to E$ uniformly on a compact set $\Omega \times [0, \infty)$ and $\Phi$ and $E$ are continuous on $\Omega \times [0, \infty)$. So $\Phi_n(E_n) \to \Phi(E)$ uniformaly on $\Omega \times [0, \infty)$. Then, for $\mu$-a.e. $t \in \Omega$, $\Phi(E(t, u)) = \lim_{n \to \infty} \Phi_n(E_n(t, u))$ is convex of $u$ on $[0, \infty)$,
$$\Phi(E(t, 0)) = \lim_{u \to 0^+} \Phi(E(t, u)) = \lim_{n \to \infty} \lim_{u \to 0^+} \Phi_n(E_n(t, u)) = 0,$$
$$\lim_{u \to \infty} \Phi(E(t, u)) = \lim_{n \to \infty} \lim_{u \to \infty} \Phi_n(E_n(t, u)) = \infty$$
and for each $u \in [0, \infty)$, $\Phi(E(t, u))$ is a $\mu$-measurable function of $t$ on $\Omega$. ∎

**Proposition 3.11.** Let $(\Omega, \Sigma, \mu)$ be a measure space and $\Phi_1, \Phi_2: \Omega \times [0, \infty) \to \mathbb{R}$ are $E$-strong Young functions with a map $E: \Omega \times [0, \infty) \to \Omega \times [0, \infty)$. Then both $\Phi_1 + \Phi_2$ and $c\Phi_1, c \geq 0$ are $E$-strong Young functions with $E$.

**Proposition 3.12.** Let $(\Omega, \Sigma, \mu)$ be a measure space with $\Omega \times [0, \infty)$ compact. If $\Phi_n: \Omega \times [0, \infty) \to \mathbb{R}$ is a sequence of continuous $E$-strong Young functions with a map $E: \Omega \times [0, \infty) \to \Omega \times [0, \infty)$ and $\Phi_n$ converges uniformaly to a continuous map $\Phi: \Omega \times [0, \infty) \to \mathbb{R}$, then $\Phi$ is an $E$-strong Young function with $E$.

**Proof.** Since $\Phi_n$ are continuous $E$-strong Young functions with $E$ such that $\Phi_n \to \Phi$ uniformly on a compact set $\Omega \times [0, \infty)$ and $\Phi$ is continuous on $\Omega \times [0, \infty)$, then $\Phi_n(E) \to \Phi(E)$ uniformaly on $\Omega \times [0, \infty)$. So, for $\mu$-a.e. $t \in \Omega$, $\Phi(E(t, u)) = \lim_{n \to \infty} \Phi_n(E(t, u))$ is convex continuous of $u$ on $[0, \infty)$,

$$\Phi(E(t, 0)) = \lim_{n \to \infty} \Phi_n(E(t, 0)) = 0 \Leftrightarrow u = 0,$$

$$\lim_{u \to \infty} \Phi(E(t, u)) = \lim_{n \to \infty} \lim_{u \to \infty} \Phi_n(E(t, u)) = \infty$$

and for each $u \in [0, \infty)$, $\Phi(E(t, u))$ is a $\mu$-measurable function of $t$ on $\Omega$. ∎

**Proposition 3.13.** Let $(\Omega, \Sigma, \mu)$ be a measure space and $\Phi: \Omega \times [0, \infty) \to \mathbb{R}$ be a linear $E$-strong Young function with maps $E_1, E_2: \Omega \times [0, \infty) \to \Omega \times [0, \infty)$. Then $\Phi$ is $E$-strong Young function with $E_1 + E_2$ and $cE_1$.

**Proposition 3.14.** Let $(\Omega, \Sigma, \mu)$ be a measure space with $\Omega \times [0, \infty)$ compact. If $\Phi: \Omega \times [0, \infty) \to \mathbb{R}$ be a continuous $E$-strong Young function with a sequence of maps $E_n: \Omega \times [0, \infty) \to \Omega \times [0, \infty)$ that converges uniformaly to a map $E: \Omega \times [0, \infty) \to \Omega \times [0, \infty)$. Then $\Phi$ is an $E$-strong Young function with $E$.

**Proof.** Because $\Phi$ is a continuous $E$-strong Young function with $E_n$ such that $E_n \to E$ uniformaly on a compact set $\Omega \times [0, \infty)$ and $E$ is continuous on $\Omega \times [0, \infty)$, then $\Phi(E_n) \to \Phi(E)$ uniformaly on $\Omega \times [0, \infty)$. So, for $\mu$-a.e. $t \in \Omega$, $\Phi(E(t, u)) = \lim_{n \to \infty} \Phi(E_n(t, u))$ is convex continuous of $u$ on $[0, \infty)$,

$$\Phi(E(t, 0)) = \lim_{n \to \infty} \Phi(E_n(t, 0)) = 0 \Leftrightarrow u = 0,$$
$$\lim_{u \to \infty} \Phi(E(t, u)) = \lim_{n \to \infty} \lim_{u \to \infty} \Phi(E_n(t, u)) = \infty$$

and for each $u \in [0, \infty)$, $\Phi(E(t, u))$ is a $\mu$-measurable function of $t$ on $\Omega$. ∎

**Proposition 3.15.** Let $(\Omega, \Sigma, \mu)$ be a measure space with $\Omega \times [0, \infty)$ compact. If $\Phi_n: \Omega \times [0, \infty) \to \mathbb{R}$ is a sequence of continuous $E$-strong Young functions with a sequence of continuous maps $E_n: \Omega \times [0, \infty) \to \Omega \times [0, \infty)$ such that $\Phi_n$ converges uniformaly to a continuous function $\Phi: \Omega \times [0, \infty) \to \mathbb{R}$ and $E_n$ converges uniformaly to a continuous map $E: \Omega \times [0, \infty) \to \Omega \times [0, \infty)$, then $\Phi$ is an $E$-strong Young function with $E$.

**Proof.** We have $\Phi_n$ are continuous $E$-strong Young functions with continuous maps $E_n$ such that $\Phi_n \to \Phi$ and $E_n \to E$ uniformly on a compact set $\Omega \times [0, \infty)$ and $\Phi$ and $E$ are continuous on $\Omega \times [0, \infty)$. So, $\Phi_n(E_n) \to \Phi(E)$ uniformaly on $\Omega \times [0, \infty)$. Then, for $\mu$-a.e. $t \in \Omega$, $\Phi(E(t, u)) = \lim_{n \to \infty} \Phi_n(E_n(t, u))$ is convex continuous of $u$ on $[0, \infty)$,

$$\Phi(E(t, 0)) = \lim_{n \to \infty} \Phi_n(E_n(t, 0)) = 0 \Leftrightarrow u = 0,$$

$$\lim_{u\to\infty} \Phi\big(E(t,u)\big) = \lim_{n\to\infty}\lim_{u\to\infty} \Phi_n\big(E_n(t,u)\big) = \infty$$

and for each $u \in [0,\infty)$, $\Phi\big(E(t,u)\big)$ is a $\mu$-measurable function of $t$ on $\Omega$. ∎

**Proposition 3.16.** Let $(\Omega, \Sigma, \mu)$ be a measure space and $\Phi_1, \Phi_2 \colon \Omega \times [0,\infty) \to \mathbb{R}$ are $E$-Orlicz functions with a map $E \colon \Omega \times [0,\infty) \to \Omega \times [0,\infty)$. Then $\Phi_1 + \Phi_2$ and $c\Phi_1, c \geq 0$ are $E$-Orlicz functions with $E$.

**Proposition 3.17.** Let $(\Omega, \Sigma, \mu)$ be a measure space with $\Omega \times [0,\infty)$ compact. If $\Phi_n \colon \Omega \times [0,\infty) \to \mathbb{R}$ is a sequence of continuous $E$-Orlicz functions with a map $E \colon \Omega \times [0,\infty) \to \Omega \times [0,\infty)$ and $\Phi_n$ converges uniformaly to a continuous map $\Phi \colon \Omega \times [0,\infty) \to \mathbb{R}$, then $\Phi$ is an $E$-Orlicz function with $E$.

**Proof.** Since $\Phi_n$ are continuous $E$-Orlicz functions with $E$ such that $\Phi_n \to \Phi$ uniformly on a compact set $\Omega \times [0,\infty)$ and $\Phi$ is continuous on $\Omega \times [0,\infty)$, then $\Phi_n(E) \to \Phi(E)$ uniformaly on $\Omega \times [0,\infty)$. So, for $\mu$-a.e. $t \in \Omega$, $\Phi\big(E(t,u)\big) = \lim_{n\to\infty} \Phi_n\big(E(t,u)\big)$ is convex of $u$ on $[0,\infty)$,

$$\Phi\big(E(t,0)\big) = \lim_{n\to\infty} \Phi_n\big(E(t,0)\big) = 0,$$

$$\lim_{u\to\infty} \Phi\big(E(t,u)\big) = \lim_{n\to\infty}\lim_{u\to\infty} \Phi_n\big(E(t,u)\big) = \infty,$$

$0 < \Phi\big(E(t,u)\big) < \infty$ for any $u \in (0,\infty)$ and $\Phi\big(E(t,u)\big)$ is left continuous at

$$U_\Phi = \sup\{u > 0 \colon \Phi\big(E(t,u)\big) < +\infty\}.$$

Furthermore, for each $u \in [0,\infty)$, $\Phi\big(E(t,u)\big)$ is a $\mu$-measurable function of $t$ on $\Omega$. ∎

**Proposition 3.18.** Let $(\Omega, \Sigma, \mu)$ be a measure space and $\Phi \colon \Omega \times [0,\infty) \to \mathbb{R}$ be a linear $E$-Orlicz function with maps $E_1, E_2 \colon \Omega \times [0,\infty) \to \Omega \times [0,\infty)$, then $\Phi$ is $E$-Orlicz function with $E_1 + E_2$ and $cE_1, c \geq 0$.

**Proposition 3.19.** Let $(\Omega, \Sigma, \mu)$ be a measure space with $\Omega \times [0,\infty)$ compact. If $\Phi \colon \Omega \times [0,\infty) \to \mathbb{R}$ is a continuous $E$-Orlicz function with a sequence of maps $E_n \colon \Omega \times [0,\infty) \to \Omega \times [0,\infty)$ that converges uniformaly to a map $E \colon \Omega \times [0,\infty) \to \Omega \times [0,\infty)$, then $\Phi$ is an $E$-Orlicz function with $E$.

**Proof.** Because $\Phi$ is a continuous $E$-Orlicz function with $E_n$ such that $E_n \to E$ uniformaly on a compact set $\Omega \times [0,\infty)$ and $E$ is continuous on $\Omega \times [0,\infty)$, then $\Phi(E_n) \to \Phi(E)$ uniformaly on $\Omega \times [0,\infty)$. So, for $\mu$-a.e. $t \in \Omega$, $\Phi\big(E(t,u)\big) = \lim_{n\to\infty} \Phi\big(E_n(t,u)\big)$ is convex of $u$ on $[0,\infty)$,

$$\Phi\big(E(t,0)\big) = \lim_{n\to\infty} \Phi\big(E_n(t,0)\big) = 0,$$

$$\lim_{u\to\infty} \Phi\big(E(t,u)\big) = \lim_{n\to\infty}\lim_{u\to\infty} \Phi\big(E_n(t,u)\big) = \infty,$$

$0 < \Phi\big(E(t,u)\big) < \infty$ for any $u \in (0,\infty)$ and $\Phi\big(E(t,u)\big)$ is left continuous at

$$U_\Phi = \sup\{u > 0 \colon \Phi\big(E(t,u)\big) < +\infty\}.$$

Moreover, for each $u \in [0,\infty)$, $\Phi\big(E(t,u)\big)$ is a $\mu$-measurable function of $t$ on $\Omega$. ∎

**Proposition 3.20.** Let $(\Omega, \Sigma, \mu)$ be a measure space with $\Omega \times [0, \infty)$ compact. If $\Phi_n: \Omega \times [0, \infty) \to \mathbb{R}$ is a sequence of continuous $E$-Orlicz functions with a sequence of continuous maps $E_n: \Omega \times [0, \infty) \to \Omega \times [0, \infty)$ such that $\Phi_n$ converges uniformaly to a continuous function $\Phi: \Omega \times [0, \infty) \to \mathbb{R}$ and $E_n$ converges uniformaly to a continuous map $E: \Omega \times [0, \infty) \to \Omega \times [0, \infty)$, then $\Phi$ is an $E$-Orlicz function with $E$.

**Proof.** We have $\Phi_n$ are continuous $E$-Orlicz functions with continuous maps $E_n$ such that $\Phi_n \to \Phi$ and $E_n \to E$ uniformly on a compact set $\Omega \times [0, \infty)$ and $\Phi$ and $E$ are continuous on $\Omega \times [0, \infty)$. So $\Phi_n(E_n) \to \Phi(E)$ uniformaly on $\Omega \times [0, \infty)$. Then, for $\mu$-a.e. $t \in \Omega, \Phi(E(t,u)) = \lim_{n \to \infty} \Phi_n(E_n(t,u))$ is convex of $u$ on $[0, \infty)$,

$$\Phi(E(t,0)) = \lim_{n \to \infty} \Phi_n(E_n(t,0)) = 0,$$

$$\lim_{u \to \infty} \Phi(E(t,u)) = \lim_{n \to \infty} \lim_{u \to \infty} \Phi_n(E_n(t,u)) = \infty,$$

$0 < \Phi(E(t,u)) < \infty$ for any $u \in (0, \infty)$ and $\Phi(E)$ is left continuous at

$$U_\Phi = \sup\{u > 0 : \Phi(E(t,u)) < +\infty\}.$$

Moreover, for each $u \in [0, \infty), \Phi(E(t,u))$ is a $\mu$-measurable function of $t$ on $\Omega$. ∎

## 4. Relationships between $E$-convex functions $\Phi$

In this section, we generalize the theorems in [8] to the $E$-ones and follow the same way to prove them to consider the relationship between $E$-$N$-function function, $E$-Young function, $E$-strong Young function and $E$-Orlicz function.

**Theorem 4.1.** If $\Phi$ is an $E$-$N$-function, then $\Phi$ is an $E$-strong Young function.

**Proof.** Let $\Phi: \Omega \times [0, \infty) \to \mathbb{R}$ be an $E$-$N$-function with a map $E: \Omega \times [0, \infty) \to \Omega \times [0, \infty)$. Then, for $\mu$-a.e. $t \in \Omega, \Phi(E(t,u))$ is convex continuous of $u$ on $[0, \infty)$ satisfying

$$\forall \varepsilon > 0, \exists \delta > 0, 0 < u < \delta \Longrightarrow \left|\frac{\Phi(E(t,u))}{u}\right| < \varepsilon$$

because $\lim_{u \to 0^+} \frac{\Phi(E(t,u))}{u} = 0$. Letting $\delta < 1$, we have

$$0 \leq |\Phi(E(t,u))| < \left|\frac{\Phi(E(t,u))}{\delta}\right| < \left|\frac{\Phi(E(t,u))}{u}\right| < \varepsilon.$$

By the squeeze theorem for functions, we get $\Phi(E(t,0)) = 0 \Leftrightarrow u = 0$ because $\Phi$ is continuous at $u = 0$ and $\Phi(E(t,u)) > 0$ for any $u \in (0, \infty)$. Also,

$$\forall M \in \mathbb{R}, \exists u_M > 0, u > u_M \Longrightarrow \frac{\Phi(E(t,u))}{u} > M$$

because $\lim_{u \to \infty} \frac{\Phi(E(t,u))}{u} = \infty$. Taking $u_M > 1$, we have that

$$\Phi(E(t,u)) > Mu > Mu_M > M.$$

That is, $\lim_{u\to\infty} \Phi(E(t,u)) = \infty$. Furthermore, for each $u \in [0,\infty)$, $\Phi(E(t,u))$ is a $\mu$-measurable function of $t$ on $\Omega$ which completes the proof. ∎

**Remark 4.1:** An $E$-strong Young function may not be $E$-$N$-function.

**Example 4.1:** Consider the function $\Phi: \mathbb{R} \times [0,\infty) \to \mathbb{R}$ defined as $\Phi(t,u) = e^{u^t} - 1$ with the map $E: \mathbb{R} \times [0,\infty) \to \mathbb{R} \times [0,\infty)$ defined as $E(t,u) = (1,u)$. Then $\Phi$ is $E$-strong Young function but is not $E$-$N$-function because for $\mu$-a.e. $t \in \mathbb{R}$, $\lim_{u\to 0} \frac{e^u - 1}{u} = 1 \neq 0$.

**Theorem 4.2.** If $\Phi$ is an $E$-strong Young function, then $\Phi$ is an $E$-Orlicz function.

**Proof.** Let $\Phi: \Omega \times [0,\infty) \to \mathbb{R}$ be an $E$-strong Young function with a map $E: \Omega \times [0,\infty) \to \Omega \times [0,\infty)$. Then for $\mu$-a.e. $t \in \Omega$, $\Phi(E(t,u))$ is continuous convex of $u$ on $[0,\infty)$ satisfying $\Phi(E(t,0)) = 0$, $\Phi(E(t,u)) > 0$ for any $u \in (0,\infty)$ because $\Phi(E(t,0)) = 0 \Leftrightarrow u = 0$ and $\lim_{u\to\infty} \Phi(E(t,u)) = \infty$ and $\Phi(E(t,u))$ is left continuous at $U_\Phi = +\infty$ because $\lim_{u\to\infty} \Phi(E(t,u)) = \infty$. Moreover, for each $u \in [0,\infty)$, $\Phi(E(t,u))$ is a $\mu$-measurable function of $t$ on $\Omega$. Hence, $\Phi$ is an $E$-Orlicz function. ∎

**Remark 4.2:** An $E$-Orlicz function may not be $E$-strong Young function.

**Example 4.2:** Consider the function $\Phi: \mathbb{R} \times [0,\infty) \to \mathbb{R}$ defined as

$$\Phi(t,u) = \begin{cases} u - |t|, & 0 \leq u < 1 \\ u + |t| - 2, & 1 \leq u \end{cases}$$

with the map $E: \mathbb{R} \times [0,\infty) \to \mathbb{R} \times [0,\infty)$ defined as $E(t,u) = (u,u)$. Then $\Phi$ is $E$-Orlicz function but not $E$-strong Young function because that for $\mu$-a.e. $t \in \Omega$, $\Phi(E(t,1)) = 0$.

**Theorem 4.3.** If $\Phi$ is an $E$-Orlicz function, then $\Phi$ is an $E$-Young function.

**Proof.** Let $\Phi: \Omega \times [0,\infty) \to \mathbb{R}$ be an $E$-Orlicz function with a map $E: \Omega \times [0,\infty) \to \Omega \times [0,\infty)$. Then for $\mu$-a.e. $t \in \Omega$, $\Phi(E(t,u))$ is convex of $u$ on $[0,\infty)$ satisfying $\Phi(E(t,0)) = 0$,

$\lim_{u\to\infty} \Phi(E(t,u)) = \infty$, $0 < \Phi(E(t,u)) < \infty$ for any $u \in (0,\infty)$ and $\Phi(E(t,u))$ is left continuous at $U_\Phi$. We only need to show that $\lim_{u\to 0^+} \Phi(E(t,u)) = 0$. In other words, we need to prove that

$$\forall \varepsilon > 0, \exists \delta_\varepsilon > 0, 0 < u < \delta_\varepsilon \Rightarrow 0 \leq \Phi(E(t,u)) < \varepsilon.$$

For arbitrary $\varepsilon > 0$, consider $a_\Phi = \inf\{u > 0: \Phi(E(t,u)) > 0\}$. If $a_\Phi > 0$, then $\Phi(E(t,u)) = 0$ for all $u \in (0, a_\Phi)$. Taking $\delta_\varepsilon = a_\Phi > 0$, then $\Phi(E(t,u)) = 0 < \varepsilon$ for all $0 < u < \delta_\varepsilon$. So, $\lim_{u\to 0^+} \Phi(E(t,u)) = 0$ If $a_\Phi = 0$, then $\Phi(E(t,u)) > 0$ for all $u > 0$ and there exists $0 < u_0 < \infty$ such that $0 < \Phi(E(t,u_0)) < \infty$. That is, for all $\varepsilon > 0$, $\exists u_\varepsilon \in (0,\infty)$ such that $0 < \Phi(E(t,u_\varepsilon)) < \infty$. If $\Phi(E(t,u_0)) < \varepsilon$, then $\Phi(E(t,u_\varepsilon)) < \infty$ for $u_\varepsilon = u_0$ and if $\Phi(E(t,u_0)) \geq \varepsilon$, then for $u_\varepsilon = \alpha u_0$ where $0 \leq \alpha = \frac{\varepsilon}{2\Phi(E(t,u_0))} < 1$ that

$$\Phi(E(t,u_\varepsilon)) = \Phi(E(t,\alpha u_0)) \leq \alpha\Phi(E(t,u_0)) \leq \frac{\varepsilon}{2} < \varepsilon$$

because $\Phi$ is $E$-convex of $u$ on $[0,\infty)$. Taking $\delta_\varepsilon = u_\varepsilon > 0$, we get

$$0 \leq \Phi(E(t,u)) \leq \Phi(E(t,\delta_\varepsilon)) = \Phi(E(t,u_\varepsilon)) < \varepsilon, 0 < u < \delta_\varepsilon$$

because $\Phi(E(t,u))$ is increasing of $u$ on $[0,\infty)$. Furthermore, for each $u \in [0,\infty)$, $\Phi(E(t,u))$ is a $\mu$-measurable function of $t$ on $\Omega$. Hence, $\Phi$ is an $E$-Young function. ∎

**Remark 4.3:** An $E$-Young function may not be $E$-Orlicz function.

**Example 4.3:** Consider the function $\Phi: [0,\infty) \times [0,\infty) \to \mathbb{R}$ defined as

$$\Phi(t,u) = \begin{cases} -\log(u + |t|^{1/p} + 1), & 0 \leq u < 1 \\ +\infty, & 1 \leq u \end{cases}$$

with the map $E: [0,\infty) \times [0,\infty) \to [0,\infty) \times [0,\infty)$ defined as $E(t,u) = (u^p, u), p \geq 1$.

Then $\Phi$ is $E$-Young function but not $E$-Orlicz function because $\Phi(E(t,u))$ is not left continuous at $U_\Phi = 1$ where $\lim_{u \to 1} \Phi(E(t,u)) = -\log(3) \neq +\infty = \Phi(E(t,1))$.

**Corollary 4.4.** $E$-$N$-function $\Rightarrow$ $E$-strong Young function $\Rightarrow$ $E$-Orlicz function $\Rightarrow$ $E$-Young function.

**Corollary 4.5.** $E$-$N$-function $\not\Leftarrow$ $E$-strong Young function $\not\Leftarrow$ $E$-Orlicz function $\not\Leftarrow$ $E$-Young function.

## 5. Applications and Results

We can define the $E$-Orlicz spaces by generalize the concept of the Orlicz spaces using the concept of the $E$-convex functions $\Phi$ as follows: Let $(\Omega, \Sigma, \mu)$ be a measure space and $\Phi: \Omega \times [0,\infty) \to \mathbb{R}$ be a $E$-convex function with a map $E: \Omega \times [0,\infty) \to \Omega \times [0,\infty)$. The $E$-Orlicz space $EL_{\Phi(E)}(\Omega, \Sigma, \mu)$ generated by $\Phi(E)$ is defined by

$$EL_{\Phi(E)}(\Omega, \Sigma, \mu) = \left\{ f \in X_\Omega : \int_\Omega \Phi(E(t, \|f(t)\|_{BS})) d\mu < \infty \right\}$$

where $X_\Omega$ is the set of all $\mu$-measurable functions from $\Omega$ to $X$ and $(X, \|\cdot\|_X)$ is a Banach space, equipped with the $E$-Luxemburg norm

$$\|f\|_{\Phi(E)} = \inf\left\{ \lambda > 0: \int_\Omega \Phi\left(E\left(t, \frac{\|f(t)\|_{BS}}{\lambda}\right)\right) d\mu \leq 1 \right\}, f \in EL_{\Phi(E)}(\Omega, \Sigma, \mu).$$

And the $E$-Orlicz-Sobolev space $EW^k L_{\Phi(E)}(\Omega, \Sigma, \mu)$ generated by $\Phi(E)$ is the space of all functions $f \in X_\Omega$ such that $f$ and its distributional (weak partial) derivatives $D^\alpha f$ up to order $k$ are in the $E$-Orlicz space $EL_{\Phi(E)}(\Omega, \Sigma, \mu)$. That is,

$$EW^k L_{\Phi(E)}(\Omega, \Sigma, \mu) = \{ f \in EL_{\Phi(E)}(\Omega, \Sigma, \mu) : D^\alpha f \in EL_{\Phi(E)}(\Omega, \Sigma, \mu), \forall |\alpha| \leq k \}$$

equipped with the norm

$$\|f\|_{k,\Phi(E)} = \sum_{|\alpha|\leq k} \|D^\alpha f\|_{\Phi(E)}, f \in EW^k L_{\Phi(E)}(\Omega, \Sigma, \mu).$$

**Example 5.1:** From example 2.2.1, we have that $\Phi: \mathbb{R} \times [0, \infty) \to \mathbb{R}, \Phi(t, u) = e^{t+u} - 1$ is $E$-Young function with $E: \mathbb{R} \times [0, \infty) \to \mathbb{R} \times [0, \infty), E(t, u) = (u, u)$. Then the $E$-Orlicz space $EL_{\Phi(E)}(\Omega, \Sigma, \mu)$ generated by $\Phi(E(t, u)) = e^{2u} - 1$ is

$$EL_{e^{2u}-1}(\Omega, \Sigma, \mu) = \left\{ f \in X_\Omega : \int_\Omega (exp(2 \| f(t) \|_{BS}) - 1) d\mu < \infty \right\}$$

equipped with the Luxemburg norm for $f \in EL_{e^{2u}-1}(\Omega, \Sigma, \mu)$,

$$\|f\|_{e^{2u}-1} = \inf\left\{ \lambda > 0 : \int_\Omega \left( exp\left( \frac{2 \| f(t) \|_{BS}}{\lambda} \right) - 1 \right) d\mu \leq 1 \right\}.$$

Moreover, we can consider the $E$-Young function $\Phi_p(E(t,u)) = u^p, p \geq 1$ as a special case to define the $EL_p(\Omega, \Sigma, \mu)$ spaces

$$EL_p(\Omega, \Sigma, \mu) = \left\{ f \in X_\Omega : \int_\Omega |f|^p d\mu < \infty \right\}$$

equipped with the norm

$$\|f\|_p = \|f\|_{\Phi_p} = \inf\left\{ \lambda > 0 : \int_\Omega |f/\lambda|^p d\mu \leq 1 \right\} = \left( \int_\Omega |f|^p d\mu \right)^{1/p}$$

and the $E$-Sobolev spaces $EW^{k,p}(\Omega, \Sigma, \mu)$ such that

$$EW^{k,p}(\Omega, \Sigma, \mu) = \{ f \in EL_p(\Omega, \Sigma, \mu) : D^\alpha f \in EL_p(\Omega, \Sigma, \mu), \forall |\alpha| \leq k \}$$

equipped with the norm

$$\|f\|_{k,p} = \left( \sum_{|\alpha|\leq k} \|D^\alpha f\|_p \right)^{1/p}, \forall f \in EL_p(\Omega, \Sigma, \mu).$$

**Example 5.2:** Let $\Phi: \mathbb{C} \times [0, \infty) \to \mathbb{R}$ be defined as

$$\Phi(t, u) = \begin{cases} t\ln(u), & u > 1 \\ 0, & 0 \leq u \leq 1 \end{cases}$$

with $E: \mathbb{C} \times [0, \infty) \to \mathbb{C} \times [0, \infty)$ be defined as

$$E(t,u) = \begin{cases} (1, e^{u^p}), & 1 \leq p \\ (1,0), & 1 < u, p = +\infty \\ (0,0), 0 \leq u \leq 1, p = +\infty \end{cases}$$

We have for $\mu$-a.e. $t \in \mathbb{C}$, that

$$\Phi\big(E(t,u)\big) = \begin{cases} u^p, & 1 \leq p \\ +\infty, 1 < u, p = +\infty \\ 0, 0 \leq u \leq 1, p = +\infty \end{cases}$$

is $E$-Young function and the obtained spaces are $EL_p(\Omega, \Sigma, \mu), 1 \leq p \leq +\infty$.